\documentclass{amsart}
\usepackage{amsthm}
\usepackage{amsfonts}
\usepackage{amssymb}
\usepackage{amsmath}

\begin{document}
\small



\title{\textbf{A notion of $\alpha\beta$-statistical convergence of order $\gamma$ in probability}}


\author{Pratulananda Das$^{(1)}$,  Sanjoy Ghosal$^{(2)}$, Vatan Karakaya$^{(3)}$, Sumit Som$^{(4)*}$,}

\thanks{$^{*}$ This work is funded by UGC Research, HRDG, India. $^{(1)}$Department of Mathematics, Jadavpur University, Kolkata-700032, West Bengal, India.  E-mail:  pratulananda@yahoo.co.in  (P. Das).
$^{(2)}$School of Sciences, Netaji Subhas Open University, Kalyani, Nadia-741235, West Bengal, India. E-mail:  sanjoykrghosal@yahoo.co.in  (S. Ghosal).
$^{(3)}$Department of Mathematical Engineering, Yildiz Technical University, Davutpasa Campus, Esenler, 34750 Istanbul, Turkey. E-mail: vkkaya@yildiz.edu.tr  (V.Karakaya).
$^{(4)}$UGC Fellow, Department of Mathematics, Jadavpur University, Kolkata-700032, West Bengal, India. E-mail:  somkakdwip@gmail.com  (S. Som).}

\date{}

\setcounter {page}{1}

\maketitle

\begin{abstract}{A sequence of real numbers $\{x_{n}\}_{n\in \mathbb{N}}$ is said to be $\alpha \beta$-statistically convergent of order $\gamma$ (where $0<\gamma\leq 1$) to a real number $x$ \cite{a} if for every $\delta>0,$
$$\underset{n\rightarrow \infty}{\lim}\frac{1}{(\beta_{n}-\alpha_{n}+1)^\gamma}~|\{k \in [\alpha_n,\beta_n] : |x_{k}-x|\geq \delta \}|=0.$$ where $\{\alpha_{n}\}_{n\in \mathbb{N}}$ and $\{\beta_{n}\}_{n\in \mathbb{N}}$ be two sequences of  positive real numbers such that  $\{\alpha_{n}\}_{n\in \mathbb{N}}$ and $\{\beta_{n}\}_{n\in \mathbb{N}}$ are both non-decreasing, $\beta_{n}\geq \alpha_{n}$ $\forall ~n\in \mathbb{N},$ ($\beta_{n}-\alpha_{n})\rightarrow \infty$ as $n\rightarrow \infty.$ In this paper we study a related concept of convergences in which the value $|x_{k}-x|$ is replaced by $P(|X_{k}-X|\geq \varepsilon)$ and $E(|X_{k}-X|^{r})$ repectively (Where $X, X_k$ are random variables for each $k\in \mathbb{N}$, $\varepsilon>0$, $P$ denote the probability, $E$ denote the expectation) and we call them $\alpha \beta$-statistical convergence of order $\gamma$ in probability and $\alpha\beta$-statistical convergence of order $\gamma$ in $r^{\mbox{th}}$ expectation respectively. The results are applied to build the probability distribution for $\alpha\beta$-strong $p$-Ces$\grave{\mbox{a}}$ro summability of order $\gamma$ in probability and $\alpha\beta$-statistical convergence of order $\gamma$ in distribution. So our main objective is to interpret a relational behavior of above mentioned four convergences. We give a condition under which a sequence of random variables will converge to a unique limit under two different $(\alpha,\beta)$ sequences and this is also use to prove that if this condition violates then the limit value of $\alpha \beta$-statistical convergence of order $\gamma$ in probability of a sequence of random variables for two different $(\alpha,\beta)$ sequences may not be equal.}
 \end{abstract}





\section{\textbf{Introduction}}
\baselineskip .82cm \small

The idea of convergence of a real sequence has been extended to statistical convergence by Fast \cite{fa} and Steinhaus \cite{st} and later on re-introduced by Schoenberg \cite{sc} independently and is based on the notion of asymptotic density of the subset of natural numbers. However, the first idea of statistical convergence (by different name) was given by Zygmund \cite{z} in the first edition of his monograph published in Warsaw in 1935. Later on it was further investigated from the sequence space point of view and linked with summability theorem by Fridy \cite{fr}, Connor \cite{co}, \v{S}al\'{a}t \cite{s}, Das et. al. \cite{pd4}, Fridy and Orhan \cite{fo}.\\

In \cite{bdp,c} a different direction was given to the study of statistical convergence where the notion of statistical convergence of order $\gamma$ ($0<\gamma<1$) was introduced by using the notion of natural density of order $\gamma$ (where $n$ is replaced by $n^{\gamma}$ in the denominator in the definition of natural density). It was observed in \cite{bdp}, that the behavior of this new convergence was not exactly parallel to that of statistical convergence and some basic properties were obtained. More results on this convergence can be seen from \cite{se}.\\

In this context it should be noted that the history of strong $p$-Ces$\grave{\mbox{a}}$ro summability is not so clear. Connor in \cite{co}, observed that if a sequence is strongly $p$-Ces$\grave{\mbox{a}}$ro summable of order $\gamma$ (for $0<p<\infty$) to $x$, then the sequence must be statistically convergent of order $\gamma$ to the same limit. Both Fast \cite{fa} and Schoenberg \cite{sc} noted that if a bounded sequence is statistically convergent to $x$, then it is strongly Ces$\grave{\mbox{a}}$ro summable to $x.$ But in the more general case of order $\gamma$ the result may not be true, as was established in \cite{c}. In \cite{fsr}, the relation between strongly Ces$\grave{\mbox{a}}$ro summable and $N_{\theta}$-convergence was established among other things.\\

Recently the idea of statistical convergence of order $\gamma$ was further extended to $\alpha\beta$-statistical convergence of order $\gamma$ in \cite{a} as follows: Let $\alpha=\{\alpha_n\}_{n\in \mathbb{N}}$, $\beta=\{\beta_n\}_{n\in \mathbb{N}}$ be two non-decreasing sequences of positive real numbers satisfying the conditions,
 $\alpha_n \leq \beta_n ~\forall~~ n\in \mathbb{N}$ and $(\beta_n-\alpha_n)\rightarrow \infty$ as $n\rightarrow \infty$. This pair of sequence we denoted by ($\alpha,\beta$). Then a sequence $\{x_{n}\}_{n\in \mathbb{N}}$  of real numbers is said to be $\alpha\beta$-statistically convergent of order $\gamma$ (where $0<\gamma\leq 1$) to a real number $x$  if for each $\varepsilon>0,$ the set $K=\{n\in \mathbb{N}: |x_n-x|\geq \varepsilon\}$ has $\alpha\beta$-natural density zero, i.e $$\displaystyle{\lim_{n\rightarrow \infty}}\frac{1}{(\beta_n-\alpha_n+1)^{\gamma}}|\{k\in [\alpha_n,\beta_n]:|x_k-x|\geq \varepsilon\}|=0$$ and we write $S_{\alpha\beta}^{\gamma}-\lim x_{n}=x$ or $x_{n}\xrightarrow {S_{\alpha\beta}^{\gamma}}x.$ $\alpha\beta$-statistical convergence of order $\gamma$ is more general than statistical convergence of order $\gamma$, lacunary statistical convergence of order $\gamma$ and $\lambda$ statistical convergence of order $\gamma$ if we take (i)$\alpha_n=1$ and $\beta_n=n, \forall n\in \mathbb{N},$ (ii)$\alpha_r=k_{r-1}+1$ and $\beta_r=k_r, \forall r\in \mathbb{N},$ where $\{k_r\}_{r\in \mathbb{N}\cup \{0\}}$ is a lacunary sequence,(iii)$\alpha_n=n-\lambda_{n}+1$ and $\beta_n=n, \forall n\in \mathbb{N}$ respectively.\\

On the other hand in probability theory, a new type of convergence called statistical convergence in probability was introduced in \cite{g1}, as follows: Let $\{X_{n}\}_{n\in \mathbb{N}}$ be a sequence of random variables where each $X_{n}$ is defined on the same sample space $S$ (for each $n$) with respect to a given class of events $\triangle$ and a given probability function  $P:\triangle \rightarrow \mathbb{R}.$ Then the sequence $\{X_{n}\}_{n\in \mathbb{N}}$ is said to be statistically convergent in probability to a random variable $X$ (where $X:S\rightarrow \mathbb{R}$) if for any $\varepsilon,\delta >0$ $$ {\lim_{n\rightarrow \infty}}\frac{1}{n}|\{k\leq n: P(|X_{k}-X|\geq \varepsilon)\geq \delta\}|=0.$$ In this case we write $X_{n}\xrightarrow{PS} X.$  The class of all  sequences of random variables which are statistically convergent in probability is denoted by $PS$. One can also see \cite{pd6,g2,pd7} for related works.\\

In a natural way, in this paper we combine the approches of the above mentioned papers and introduce new and more general methods, namely,
$\alpha \beta$-statistical convergence of order $\gamma$ in probability, $\alpha\beta$-strong $p$-Ces$\grave{\mbox{a}}$ro summability of order $\gamma$ in probability, $\alpha\beta$-statistical convergence of order $\gamma$ in $r^{\mbox{th}}$ expectation and $\alpha\beta$-statistical convergence of order $\gamma$ in distribution. We mainly investigate their relationship and also make some observations about these classes. In the way we show that $\alpha\beta$-statistical limit of order $\gamma$ ($0<\gamma<1$) of a sequence of random variables for two different ($\alpha,\beta$) sequences may not be equal. It is important to note that the method of proofs and in particular examples are not analogous to the real case.\\

\section{\textbf{$\alpha\beta$-statistical convergence of order $\gamma$ in probability}}

 We first introduce the definition of $\alpha\beta$-statistical convergence of order $\gamma$ in probability for a sequence of random variables as follows.\\

\noindent{\textbf{Definition 2.1.}} Let $(S,\triangle,P)$ be a probability space and $\{X_{n}\}_{n\in \mathbb{N}}$ be a sequence of random variables where each $X_{n}$ is defined on the same sample space $S$ (for each $n$) with respect to a given class of events $\triangle$ and a given probability function  $P:\triangle \rightarrow \mathbb{R}.$ Then the sequence $\{X_{n}\}_{n\in \mathbb{N}}$ is said to be $\alpha\beta$-statistically convergent of order $\gamma$ (where $0<\gamma\leq 1$) in probability to a random variable $X$ (where $X:S\rightarrow \mathbb{R}$) if for any $\varepsilon,\delta >0$ $$ {\lim_{n\rightarrow \infty}}\frac{1}{(\beta_n-\alpha_n+1)^{\gamma}}|\{k\in [\alpha_n,\beta_n]: P(|X_{k}-X|\geq \varepsilon)\geq \delta\}|=0$$ or equivalantly $$ {\lim_{n\rightarrow \infty}}\frac{1}{(\beta_n-\alpha_n+1)^{\gamma}}|\{k\in [\alpha_n,\beta_n]: 1-P(|X_{k}-X|<\varepsilon)\geq \delta\}|=0.$$ In this case we write $S_{\alpha\beta}^{\gamma}-\lim P(|X_{n}-X|\geq \varepsilon)=0$  or  $S_{\alpha\beta}^{\gamma}-\lim P(|X_{n}-X|< \varepsilon)=1$ or just $X_{n}\xrightarrow{PS_{\alpha\beta}^{\gamma}} X.$  The class of all  sequences of random variables which are $\alpha\beta$-statistically convergent of order $\gamma$ in probability is denoted simply by $PS_{\alpha\beta}^{\gamma}.$\\

\noindent{\textbf{Note 2.1.}} In Definition 2.1 if we take $\alpha_{n}=1$ and $\beta_{n}=n$, then $\{X_{n}\}_{n\in \mathbb{N}}$ is said to be statistically convergent of order $\gamma$ (where $0<\gamma\leq 1$) in probability to a random variable $X$. So $\alpha\beta$-statistical convergence of order $\gamma$ in probability is a generalization of statistical convergence of order $\gamma$ in probability for a sequence of random variables.\\

 To show that this is indeed more general we will now give an example of a sequence of random variables which is $\alpha\beta$-statistically convergent of order $\frac{1}{2}$ in probability but is not statistically convergent of order $\frac{1}{2}$ in probability.\\

\noindent{\textbf{Example 2.1.}}
Let a sequence of random variables ${\{X_{n}}\}_{n\in \mathbb{N}}$ be defined by
\begin{equation*}
X_{n}\in
\begin{cases}
\{-1,1\} ~~ \mbox{with probability} ~~ \frac{1}{2},~~ \mbox{if} ~~ n=m^2~~\mbox{for some}~~m\in \mathbb{N}\\
 \{0,1\} ~~ \mbox{with probability} ~ P(X_{n}=0)=(1-\frac{1}{n})~~ \mbox{and} ~~ P(X_{n}=1)=\frac{1}{n},~~\mbox{if}~~ n\neq m^2~~\mbox{for any}~~m \in \mathbb{N}
\end{cases}
\end{equation*}
Let $0<\varepsilon,\delta<1. ~~\mbox{Then we have},$
$$P(|X_{n}-0|\geq \varepsilon)=1 ~~ \mbox{if} ~~ n=m^2~~\mbox{for some}~~m\in \mathbb{N}$$ and
$$P(|X_{n}-0|\geq \varepsilon)=\frac{1}{n},~~\mbox{if}~~n\neq m^2~~\mbox{for any}~~m \in \mathbb{N}.$$

Let $\gamma=\frac{1}{2},\alpha_{n}=((n-1)^2+1), \beta_{n}=n^2 ~\forall ~ n\in \mathbb{N}.$ Then we have the innequality, $$\frac{1}{\sqrt{2n-1}}|\{k\in [(n-1)^2+1,n^2]: P(|X_{n}-0|\geq \varepsilon)\geq \delta\}|=(\frac{1}{\sqrt{2n-1}}+\frac{d}{\sqrt{2n-1}})\rightarrow 0 ~~\mbox{as}~~n\rightarrow \infty$$
where $d$ is a finite positive integer. So $X_{n}\xrightarrow {PS_{\alpha\beta}^{\frac{1}{2}}}0$.\\

But $$\frac{\sqrt{n}-1}{\sqrt{n}}\leq \frac{1}{\sqrt{n}}|\{k\leq n: P(|X_{n}-0|\geq \varepsilon)\geq \delta\}|.$$ So the right hand side does not tend to $0$. This shows that $\{X_{n}\}_{n\in \mathbb{N}}$ is not statistically convergent of order $\frac{1}{2}$ in probability to $0$.\\

\noindent{\textbf{Theorem 2.1.}} If a sequence of constants $x_{n}\xrightarrow{S_{\alpha\beta}^{\gamma}} x$ then regarding a constant as a random variable having one point distribution at that point, we may also write $x_{n}\xrightarrow{PS_{\alpha\beta}^{\gamma}}x$.\\

\noindent{\textbf{Proof :}}  Let $\varepsilon >0$ be any arbitrarily small positive real number. Then $$\displaystyle{\lim_{n\rightarrow \infty}}~\frac{1}{(\beta_n-\alpha_n+1)^{\gamma}}|\{k\in[\alpha_n,\beta_n]:|x_{k}-x|\geq \varepsilon\}|=0.$$ Now let $\delta>0$. So the set $K_{1}=\{k \in  \mathbb{N}:P(|x_{k}-x|\geq \varepsilon)\geq \delta\}\subseteq K$ where $K=\{k\in \mathbb{N}:|x_{k}-x|\geq \varepsilon\}$. This shows that  $x_{n}\xrightarrow{PS_{\alpha\beta}^{\gamma}} x$.\\

The following example shows that in general the converse of Theorem 2.1 is not true and also shows that there is a sequence $\{X_{n}\}_{n\in \mathbb{N}}$ of random variables which is $\alpha\beta$-statistically convergent in probability to a random variable X but it is not $\alpha\beta$-statistically convergent of order $\gamma$ in probability for $0<\gamma<1.$\\

\noindent{\textbf{Example 2.2.}} Let $c$ be a rational number between $\gamma_1$ and $\gamma_2$. Let the probability density function of $X_{n}$ be given by,

\begin{equation*}
f_{n}(x)=
\begin{cases}
1  ~\mbox{where} ~ 0<x<1 ~&\\
0  ~\mbox{otherwise}, ~\mbox{if} ~ n=[m^{\frac{1}{c}}]~ \mbox{for some} ~ m\in \mathbb{N}
\end{cases}
\end{equation*}
\begin{equation*}
f_{n}(x)=
\begin{cases}
\frac{nx^{n-1}}{2^{n}} ~ \mbox{where} ~~ 0<x<2~&\\
0 ~ \mbox{otherwise, if} ~ n\neq [m^{\frac{1}{c}}]~ \mbox{for any} ~~ m\in \mathbb{N}
\end{cases}
\end{equation*}
Now let $0<\varepsilon,\delta<1.$ Then
$$P(|X_{n}-2|\geq \varepsilon)=1 ~ \mbox{if} ~ n=[m^{\frac{1}{c}}] ~ \mbox{for some} ~~  m\in \mathbb{N}.$$
$$P(|X_{n}-2|\geq \varepsilon)=(1- \frac{\varepsilon}{2})^{n} ~ \mbox{if} ~ n\neq[m^{\frac{1}{c}}] ~~ \mbox{for any} ~~  m\in \mathbb{N}.$$
Now let $\alpha_n=1, ~\beta_n=n^{2}$. Consequently we have the inequality, $$ {\lim_{n\rightarrow \infty}}\frac{n^{2c}-1}{n^{2\gamma_1}} \leq  {\lim_{n\rightarrow \infty}} \frac{1}{n^{2\gamma_1}}|{\{k\in [1,n^2]:P(|X_{k}-2|\geq \varepsilon)\geq \delta \}}|$$
and $${\lim_{n\rightarrow \infty}} \frac{1}{n^{2\gamma_2}}|\{k\in [1,n^2]:P(|X_{k}-2|\geq \varepsilon)\geq \delta \}|\leq {\lim_{n\rightarrow \infty}}(\frac{n^{2c}+1}{n^{2\gamma_2}}+\frac{d}{n^{2\gamma_2}})$$
where d is a fixed finite positive integer. This shows that $\{X_{n}\}_{n\in \mathbb{N}}$ is $\alpha\beta$-statistically convergent of order $\gamma_2$ in probability to $2$ but is not $\alpha\beta$-statistically convergent of order $\gamma_1$ in probability to $2$ whenever $\gamma_1<\gamma_2$ and this is not the usual $\alpha\beta$-statistical convergence of order $\gamma$ of real numbers. So the converse of Theorem 2.1. is not true. Also by taking $\gamma_2 =1,$ we see that $X_{n}\xrightarrow {PS_{\alpha\beta}}2$ but $\{X_{n}\}_{n\in \mathbb{N}}$ is not $\alpha\beta$-statistically convergent of order $\gamma$ in probability to $2$ for $0<\gamma<1.$\\

\noindent{\textbf{Theorem 2.2.}} [Elementary properties].

(i) If $X_{n}\xrightarrow {PS_{\alpha\beta}^{\gamma_1}} X$ and $X_{n}\xrightarrow {PS_{\alpha\beta}^{\gamma_2}} Y$ then $P\{X=Y\}=1$ for any $\gamma_1,\gamma_2$ where $0<\gamma_1,\gamma_2 \leq 1$

(ii) If $X_{n}\xrightarrow {PS_{\alpha\beta}^{\gamma_1}} X$ and $Y_{n}\xrightarrow {PS_{\alpha\beta}^{\gamma_2}} Y$ then $(cX_{n}+dY_{n})\xrightarrow {PS_{\alpha\beta}^{\max\{\gamma_1,\gamma_2\}}} (cX+dY)$   where $c,d$ are constants and $0<\gamma_1,\gamma_2 \leq 1.$

(iii) Let $0<\gamma_1\leq\gamma_2\leq 1.$ Then $PS_{\alpha\beta}^{\gamma_1}\subseteq PS_{\alpha\beta}^{\gamma_2}$ and this inclusion is strict whenever $\gamma_1<\gamma_2$.

(iv) Let $g:\mathbb{R}\rightarrow\mathbb{R}$ be a continuous function and $0<\gamma_1\leq \gamma_2\leq1$. If $X_{n}\xrightarrow {PS_{\alpha\beta}^{\gamma_1}}X$ then $g(X_{n})\xrightarrow {PS_{\alpha\beta}^{\gamma_2}}g(X).$\\



\noindent{\textbf{Proof :}} (i) Without loss of generality we assume $\gamma_2\leq \gamma_1.$ If possible let $P\{X=Y\}\neq 1.$ Then there exists two positive real numbers $\varepsilon$ and $\delta$ such that $P(|X-Y|\geq \varepsilon)=\delta>0.$ Then we have
$$\displaystyle{\lim_{n\rightarrow \infty}}\frac{\beta_n-\alpha_n+1}{(\beta_n-\alpha_n+1)^{\gamma_1}}\leq {\lim_{n\rightarrow \infty}}\frac{1}{(\beta_n-\alpha_n+1)^{\gamma_1}}|\{k\in [\alpha_n,\beta_n]:P(|X_{k}-X|\geq \frac{\varepsilon}{2})\geq \frac{\delta}{2}\}|$$
$$+{\lim_{n\rightarrow \infty}}\frac{1}{(\beta_n-\alpha_n+1)^{\gamma_2}}|\{k\in [\alpha_n,\beta_n]:P(|X_{k}-Y|\geq \frac{\varepsilon}{2})\geq \frac{\delta}{2}\}|$$ which is impossible because the left hand limit is not $0$ whereas the right hand limit is $0.$ So $P\{X=Y\}=1.$  \\

Proof of (ii) is straightforward and so is omitted.\\

(iii) The first part is obvious. The inclusion is proper as can be seen from Example 2.2.\\

Proof of (iv) is straightforward and so is omitted.\\

\noindent{\textbf{Remark 2.1.}  In Theorem 2 \cite{bdp} it was observed that $m^{\gamma_1}_0\subset m^{\gamma_2}_0 $ and this inclusion was shown to be strict for at least those $\gamma_1, \gamma_2$ for which there is a $k\in \mathbb{N}$ such that $\gamma_1< \frac{1}{k}<\gamma_2.$ But Example 2.2 shows that the inequality is strict whenever $\gamma_1<\gamma_2.$\\


\noindent{\textbf{Theorem 2.3.}} Let $0<\gamma\leq 1,$ $(\alpha,\beta)$ and $(\alpha',\beta')$ are two pairs of sequences of positive real numbers such that $[\alpha'_{n},\beta'_{n}]\subseteq [\alpha_{n},\beta_{n}]$ $\forall~~ n \in \mathbb{N}$ and $(\beta_{n}-\alpha_{n}+1)^{\gamma}\leq \varepsilon (\beta'_{n}-\alpha'_{n}+1)^{\gamma}$ for some $\varepsilon>0.$ Then we have $PS_{\alpha\beta}^{\gamma}\subseteq PS_{\alpha'\beta'}^{\gamma}.$\\

\noindent{\textbf{Proof :}} Proof is straightforward and so is omitted.\\

But if the condition of the Theorem 2.3. is violated then limit may not be unique for two different $(\alpha,\beta)$'s. We now give an example to show this.\\

\noindent{\textbf{Example 2.3.} Let $\alpha=\{(2n)!\} , \beta=\{(2n+1)!\}$ and $\alpha'=\{(2n+1)!\} ,\beta'=\{(2n+2)!\}$\\
Let us define a sequence of random variables $\{X_{n}\}_{n\in \mathbb{N}}$ by,
\begin{equation*}
X_{k}\in
\begin{cases}
\{-1,1\}~~ \mbox{with probability}~~ P(X_{k}=-1)=\frac{1}{k}, P(X_{k}=1)=(1-\frac{1}{k}),~~ \mbox{if}~~(2n)!<k<(2n+1)! ~&\\
 \{-2,2\}~~ \mbox{with probability}~~ P(X_{k}=-2)=\frac{1}{k},P(X_{n}=2)=(1-\frac{1}{k}),~~ \mbox{if}~~(2n+1)!<k<(2n+2)! ~&\\
 \{-3,3\}~~ \mbox{with probability}~~ P(X_{k}=-3)=P(X_{k}=3),~~ \mbox{if}~~k=(2n)!~~\mbox{and}~~k=(2n+1)! ~&\\
\end{cases}
\end{equation*}

Let $0<\varepsilon,\delta<1$ and $0<\gamma<1.$ Then for the sequence $(\alpha,\beta)$, $$P(|X_{k}-1|\geq \varepsilon)=\frac{1}{k}~~ \mbox{if}~~(2n)!<k<(2n+1)!$$ and $$P(|X_{k}-1|\geq \varepsilon)=1~~ \mbox{if}~~(2n+1)!<k<(2n+2)!$$ and $$P(|X_{k}-1|\geq \varepsilon)=1~~ \mbox{if}~~k=(2n)!~~\mbox{and}~~k=(2n+1)!.$$

$$\Rightarrow~\displaystyle{\lim_{n\rightarrow \infty}}\frac{1}{((2n+1)!-(2n)!+1)^{\gamma}}|\{k \in [(2n)!,(2n+1)!]:P(|X_{k}-1|\geq \varepsilon)\geq \delta\}|=0$$
So $X_{n}\xrightarrow {PS_{\alpha\beta}^{\gamma}}1.$\\
Similarly it can be shown that for the sequence $\alpha'=\{(2n+1)!\} ,\beta'=\{(2n+2)!\}$, $X_{n}\xrightarrow {PS_{\alpha'\beta'}^{\gamma}}2.$\\


\noindent{\textbf{Definition 2.2.} Let $(S,\triangle,P)$ be a probability space and $\{X_{n}\}_{n\in \mathbb{N}}$ be a sequence of random variables where each $X_{n}$ is defined on the same sample space $S$ (for each $n$) with respect to a given class of events $\triangle$ and a given probability function  $P:\triangle \rightarrow \mathbb{R}.$   A sequence of random variables ${\{X_{n}\}}_{n\in \mathbb{N}}$ is said to be $\alpha\beta$-strong $p$-Ces$\grave{\mbox{a}}$ro summable of order $\gamma$ (where $0<\gamma \leq 1$ and $p>0$ is any fixed positive real number) in probability to a random variable X  if for any $\varepsilon>0$ $${\lim_{n\rightarrow\infty}}\frac{1}{(\beta_n-\alpha_n+1)^{\gamma}}\displaystyle{\sum_{k\in[\alpha_n,\beta_n]}}\{P(|X_{k}-X|\geq \varepsilon)\}^{p}=0.$$ In this case we write $X_{n}\xrightarrow {PW_{\alpha\beta}^{\gamma,p}} X.$ The class of all sequences of random variables which are $\alpha\beta$-strong $p$-Ces$\grave{\mbox{a}}$ro summable of order $\gamma$ in probability is denoted simply by $PW_{\alpha\beta}^{\gamma,p}.$\\


\noindent{\textbf{Theorem 2.4.} (i) Let $0<\gamma_1\leq \gamma_2\leq 1.$ Then $PW_{\alpha\beta}^{\gamma_1,p}\subseteq PW_{\alpha\beta}^{\gamma_2,p}.$ This inclusion is strict whenever $\gamma_1<\gamma_2.$\\
(ii) Let $0<\gamma\leq 1$ and $0<p<q<\infty$. Then $PW_{\alpha\beta}^{\gamma,q}\subset PW_{\alpha\beta}^{\gamma,p}.$\\

\noindent{\textbf{Proof :}} (i) The first part of this theorem is straightforward and so is omitted. For the second part we will give an example to show that there is a sequence of random variables ${\{X_{n}\}}_{n\in \mathbb{N}}$ which is $\alpha\beta$-strong $p$-Ces$\grave{\mbox{a}}$ro summable of order $\gamma_2$ in probability to a random variable X but is not $\alpha\beta$-strong $p$-Ces$\grave{\mbox{a}}$ro summable of order $\gamma_1$ in probability whenever $\gamma_1<\gamma_2.$\\

Let $c$ be a rational number between $\gamma_1$ and $\gamma_2$. We consider a sequence of random variables :
\begin{equation*}
X_{n} \in
\begin{cases}
 \{-1,1\}~ \mbox{with  probability} ~ \frac{1}{2},~
\mbox{if} ~ n=[m^\frac{1}{c}] ~~ \mbox{for some} ~ m\in \mathbb{N}~&\\
 \{0,1\} ~~ \mbox{with probability} ~~ P(X_{n}=0)=1-\frac{1}{\sqrt[p]{n^2}} ~~ \mbox{and} ~~ P(X_{n}=1)=\frac{1}{\sqrt[p]{n^2}}, ~ \mbox{if} ~~ n\neq [m^\frac{1}{c}] ~~\\ \mbox{for any} ~ m\in \mathbb{N}
\end{cases}
\end{equation*}
Then we have, for $0<\varepsilon<1$
$$P(|X_{n}-0|\geq \varepsilon)=1, ~ \mbox{if} ~~ n=[m^\frac{1}{c}] ~~ \mbox{for some} ~ m\in \mathbb{N}$$ and
$$P(|X_{n}-0|\geq \varepsilon)=\frac{1}{\sqrt[p]{n^2}}, ~ \mbox{if} ~~ n\neq [m^\frac{1}{c}] ~ \mbox{for any}~ m\in \mathbb{N}.$$
Let $\alpha_n=1 ~\mbox{and}~ \beta_n=n^2$. So we have the inequality $${\lim_{n\rightarrow \infty}}\frac{n^{2c}-1}{n^{2\gamma_1}}\leq {\lim_{n\rightarrow \infty}}\frac{1}{n^{2\gamma_1}}\displaystyle{\sum_{k\in [1,n^{2}]}}{\{P(|X_{k}-0|\geq \varepsilon)}\}^{p}$$ and
$${\lim_{n\rightarrow \infty}}\frac{1}{n^{2\gamma_2}}\displaystyle{\sum_{k\in [1,n^{2}]}}{\{P(|X_{k}-0|\geq \varepsilon)}\}^{p}\leq {\lim_{n\rightarrow \infty}} [\frac{n^{2c}+1}{n^{2\gamma_2}} + \frac{1}{n^{2\gamma_2}}(\frac{1}{1^{2}}+\frac{1}{2^{2}}+...+\frac{1}{n^{4}})].$$
This shows that $X_{n}\xrightarrow {PW_{\alpha\beta}^{\gamma_2,p}} 0$ but $\{X_{n}\}_{n\in \mathbb{N}}$ is not $\alpha\beta$-strong $p$-Ces$\grave{\mbox{a}}$ro summable of order $\gamma_1$ in probability to $0$.\\

(ii) Proof is straightforward and so is omitted.\\

\noindent{\textbf{Theorem 2.5.} Let $0<\gamma_1\leq\gamma_2\leq 1.$ Then $PW_{\alpha\beta}^{\gamma_1,p}\subset PS_{\alpha\beta}^{\gamma_2}.$\\

\noindent{\textbf{Proof :}} Proof is straightforward and so is omitted.\\

\noindent{\textbf{Note 2.2.}  If a sequence of random variables ${\{X_{n}}\}_{n\in \mathbb{N}}$ is $\alpha\beta$-strong $p$-Ces$\grave{\mbox{a}}$ro summable of order $\gamma$ in probability to X then it is $\alpha\beta$-statistically convergent of order $\gamma$ in probability to X i.e $PW_{\alpha\beta}^{\gamma,p}\subset PS_{\alpha\beta}^{\gamma}.$\\

But the converse of Theorem 2.5. (or Note 2.2) is not generally true as can be seen from the following example.\\

\noindent{\textbf{Example 2.4.} Let a sequence of random variables ${\{X_{n}}\}_{n\in \mathbb{N}}$ be defined by,\\
\begin{equation*}
X_{n}\in
\begin{cases}
\{-1,1\} ~~\mbox{with probability} ~~ \frac{1}{2},~~ \mbox{if} ~~ n=m^{m} ~~ \mbox{for some}~~ m\in \mathbb{N}~&\\
\{0,1\} ~~ \mbox{with probability}~~ P(X_{n}=0)=1-\frac{1}{\sqrt[2p]{n}},~~ P(X_{n}=1)=\frac{1}{\sqrt[2p]{n}},~~ \mbox{if}~ n\neq m^{m} ~ \mbox{for any} ~ m\in \mathbb{N}
\end{cases}
\end{equation*}
Let $0<\varepsilon<1$ be given. Then
$$P(|X_{n}-0|\geq \varepsilon)=1 ~~  \mbox{if} ~~ n=m^{m} ~~\mbox{for some}~~ m\in \mathbb{N}$$ and
$$P(|X_{n}-0|\geq \varepsilon)=\frac{1}{\sqrt[2p]{n}} ~~  \mbox{if} ~~ n\neq m^{m} ~~\mbox{for any}~~ m\in \mathbb{N}.$$
Let $\alpha_n=1 ~\mbox{and}~ \beta_n=n^2$. This implies $X_{n}\xrightarrow {PS_{\alpha\beta}^{\gamma}} 0$ for each $0<\gamma\leq 1.$
Next let $H={\{n\in \mathbb{N}: ~~ n\neq m^{m} ~~ \mbox{for any}~ m\in \mathbb{N}}\}.$
$$\mbox{Then}~ \frac{1}{n^{2\gamma}}\displaystyle{\sum_{k\in [1,n^2]}}{\{P(|X_{k}-0|\geq \varepsilon)}\}^{p}=\frac{1}{n^{2\gamma}}\underset{k\in H}{\sum_{k\in [1,n^2]}}{\{P(|X_{k}-0|\geq \varepsilon)}\}^{p}+\frac{1}{n^{2\gamma}}\underset{k\notin H}{\sum_{k\in [1,n^2]}}{\{P(|X_{k}-0|\geq \varepsilon)}\}^{p}$$
$$= \frac{1}{n^{2\gamma}}\underset{k\in H}{\sum_{k\in [1,n^2]}}\frac{1}{\sqrt{k}} +\frac{1}{n^{2\gamma}}\underset{k\notin H}{\sum_{k\in [1,n^2]}}1 >\frac{1}{n^{2\gamma}}{\underset{k=1}{\sum^{n^2}}}\frac{1}{\sqrt{k}}> \frac{1}{n^{2\gamma-1}} ~~(\mbox{since}~~{\underset{k=1}{\sum^{n}}}\frac{1}{\sqrt{k}}>\sqrt{n}~~ \mbox{for}~~ n\geq 2).$$
So $X_{n}$ is not $\alpha\beta$-strong $p$-Ces$\grave{\mbox{a}}$ro summable of order $\gamma$ in probability  to $0$ for $0<\gamma\leq \frac{1}{2}$.\\

\noindent{\textbf{Theorem 2.6.} (i) For $\gamma=1,$ $PW_{\alpha\beta}^{1,p}=PS_{\alpha\beta}.$\\
(ii) Let $g:\mathbb{R}\rightarrow\mathbb{R}$ be a continuous function and $0<\gamma_1\leq \gamma_2\leq1$. If $X_{n}\xrightarrow {PW_{\alpha\beta}^{\gamma_1}}X$ then $g(X_{n})\xrightarrow {PW_{\alpha\beta}^{\gamma_2}}g(X).$\\


\noindent{\textbf{Proof :}} For (i) and (ii) Proof is straightforward and so is omitted.\\

\noindent{\textbf{Theorem 2.7.} Let $\{\alpha_n\}_{n\in \mathbb{N}}$, $\{\beta_n\}_{n\in \mathbb{N}}$ be two non-decreasing sequences of positive numbers such that $\alpha_n\leq \beta_n\leq \alpha_{n+1}\leq \beta_{n+1}$ and $0<\gamma_1\leq \gamma_2\leq 1$. Then $PS^{\gamma_1}\subset PS_{\alpha\beta}^{\gamma_2}$ iff $\liminf (\frac{\beta_n}{\alpha_n})>1.$ (Such a pair of sequence exists. Take $\alpha_n=n!$ and $\beta_n=(n+1)!$).\\

\noindent{\textbf{Proof :}} First of all suppose that $\liminf (\frac{\beta_n}{\alpha_n})>1$ and let $X_{n}\xrightarrow {PS^{\gamma_1}}X$. As $\liminf (\frac{\beta_n}{\alpha_n})>1$, for each $\delta>0$ we can find sufficiently large $r$ such that $$\frac{\beta_r}{\alpha_r}\geq (1+\delta)$$ $$\Rightarrow (\frac{\beta_r-\alpha_r}{\beta_r})^{\gamma_1}\geq (\frac{\delta}{1+\delta})^{\gamma_1}.$$ Now for each $\varepsilon,\delta>0$ we have $$\frac{1}{[\beta_n]^{\gamma_1}}|\{k\leq [\beta_n]:P(|X_{k}-X|\geq \varepsilon)\geq \delta\}|= \frac{1}{[\beta_n]^{\gamma_1}}|\{k\leq \beta_n:P(|X_{k}-X|\geq \varepsilon)\geq \delta\}|$$ $$\geq\frac{1}{\beta^{\gamma_1}_n}|\{k\leq \beta_n:P(|X_{k}-X|\geq \varepsilon)\geq \delta\}|$$ $$\geq (\frac{\delta}{1+\delta})^{\gamma_1}.\frac{1}{(\beta_n-\alpha_n+1)^{\gamma_2}}|\{k\in [\alpha_n,\beta_n]:P(|X_{k}-X|\geq \varepsilon)\geq \delta\}|.$$ Hence the result follows.\\

Now if possible, suppose that $\liminf (\frac{\beta_n}{\alpha_n})=1.$ So for each $j\in \mathbb{N}$ we can choose a subsequence such that $\frac{\beta_{r(j)}}{\alpha_{r(j)}}<1+\frac{1}{j}~~ \mbox{and}~~ \frac{\beta_{r(j)-1}}{\beta_{r(j-1)}}\geq j$. Let $I_{r(j)}=[\alpha_{r(j)},\beta_{r(j)}].$\\

We define a sequence of random variables by
\begin{equation*}
X_{n}\in
\begin{cases}
 \{-1,1\}~\mbox{with probability}~ \frac{1}{2}~~\mbox{if}~ n\in I_{r(j)}~~\mbox{where}~ j\in \mathbb{N}~&\\
 \{0,1\}~\mbox{with probability}~ P(X_{n}=0)=1-\frac{1}{n^{2}}~\mbox{and}~P(X_{n}=1)=\frac{1}{n^{2}}~~\mbox{if}~ n\notin I_{r(j)}~~\mbox{for any}~ j\in \mathbb{N}
\end{cases}
\end{equation*}
Let $0<\varepsilon,\delta<1.$ Now
$$P(|X_{n}-0|\geq \varepsilon)=1 ~~\mbox{if}~ n\in I_{r(j)}~~\mbox{where}~ j\in \mathbb{N},$$ and
$$P(|X_{n}-0|\geq \varepsilon)=\frac{1}{n^{2}} ~~\mbox{if}~ n\notin I_{r(j)}~~\mbox{for any} j\in \mathbb{N}.$$\\

Now $\frac{1}{(\beta_{r(j)}-\alpha_{r(j)}+1)^{\gamma}}|\{k\in [\alpha_{r(j)},\beta_{r(j)}]: P(|X_{n}-0|\geq \varepsilon)\geq \delta\}|= \frac{(\beta_{r(j)}-\alpha_{r(j)}+1)}{(\beta_{r(j)}-\alpha_{r(j)}+1)^{\gamma}}\rightarrow \infty~~ \mbox{as}~~ j\rightarrow \infty.$ But as $\frac{1}{(\beta_{r(j)}-\alpha_{r(j)}+1)^{\gamma}}|\{k\in [\alpha_{r(j)},\beta_{r(j)}]: P(|X_{n}-0|\geq \varepsilon)\geq \delta\}|$ is a subsequence of the sequence $\frac{1}{(\beta_{r}-\alpha_{r}+1)^{\gamma}}|\{k\in [\alpha_{r},\beta_{r}]: P(|X_{n}-0|\geq \varepsilon)\geq \delta\}|,$ this shows that $X_{n}$ is not $\alpha\beta$-statistically convergent of order $\gamma$(where $0<\gamma\leq 1$) in probability to $0$.\\

Finally let $\gamma=1.$  If we take t sufficiently large such that $\alpha_{r(j)}<t\leq \beta_{r(j)}$ then we observe that,\\
$\frac{1}{t}\displaystyle{\sum_{k=1}^{t}}P(|X_{k}-0|\geq \varepsilon)\leq \frac{\beta_{r(j-1)}+\beta_{r(j)}-\alpha_{r(j)}}{\alpha_{r(j)}}+\frac{1}{t}\{1+ \frac{1}{2^{2}}+...+\frac{1}{t^{2}}\}\leq\frac{\beta_{r(j-1)}}{\beta_{r(j)-1}}+ \frac{\beta_{r(j)}-\alpha_{r(j)}}{\alpha_{r(j)}}+\frac{1}{t}\{1+ \frac{1}{2^{2}}+...+\frac{1}{t^{2}}\}\leq \frac{2}{j}+\frac{1}{t}\{1+ \frac{1}{2^{2}}+...+\frac{1}{t^{2}}\}\rightarrow 0~ \mbox{if}~ j,t\rightarrow \infty.$\\

This shows that $X_{n}\xrightarrow {PS}0.$ But this is a contradction as $PS^{\gamma_1}\subset PS_{\alpha\beta}^{\gamma_2}~~\mbox{where}(0<\gamma_1\leq \gamma_2\leq 1).$  We conclude that $\liminf (\frac{\beta_n}{\alpha_n})~\mbox{must be}>1.$\\



\noindent{\textbf{Theorem 2.8.}} Let $\alpha=\{\alpha_n\}_{n\in \mathbb{N}}$, $\beta=\{\beta_n\}_{n\in \mathbb{N}}$ be two non-decreasing sequences of positive numbers. Let $X_{n}\xrightarrow {PS^{\gamma_1}} X$ and $X_{n}\xrightarrow {PS_{\alpha\beta}^{\gamma_2}}Y$ for $0<\gamma_2\leq \gamma_1\leq 1$. If $\liminf(\frac{\beta_{n}}{\alpha_{n}})>1$ then $P\{X=Y\}=1.$\\

\noindent{\textbf{Proof :}} Let $\varepsilon> 0$ be any small positive real number and if possible let $P(|X-Y|\geq \varepsilon)=\delta>0.$ ~~ Now we have the inequality $P(|X-Y|\geq \varepsilon)\leq \{P(|X_{n}-X|\geq \frac{\varepsilon}{2})\}+ \{P(|X_{n}-Y|\geq \frac{\varepsilon}{2})\}.$\\

So ${\{k\in [\alpha_{n},\beta_{n}]:P(|X-Y|\geq \varepsilon)\geq \delta}\} \subseteq {\{k\in [\alpha_{n},\beta_{n}] : P(|X_{k}-X|\geq \frac{\varepsilon}{2})\geq \frac{\delta}{2}}\}$
$$ \bigcup {\{k\in [\alpha_{n},\beta_{n}]: P(|X_{k}-Y|\geq \frac{\varepsilon}{2})\geq \frac{\delta}{2}}\}$$\\

$\Rightarrow |{\{k\in [\alpha_{n},\beta_{n}]:P(|X-Y|\geq \varepsilon)\geq \delta}\}|\leq |{\{k\in [\alpha_{n},\beta_{n}]: P(|X_{k}-X|\geq \frac{\varepsilon}{2})\geq \frac{\delta}{2}}\}|$\\
$$+|{\{k\in [\alpha_{n},\beta_{n}]: P(|X_{k}-Y|\geq \frac{\varepsilon}{2})\geq\frac{\delta}{2}}\}|$$\\

$\Rightarrow |{\{k\in [\alpha_{n},\beta_{n}]:P(|X-Y|\geq \varepsilon)\geq \delta}\}|\leq |{\{k\leq [\beta_{n}]: P(|X_{k}-X|\geq \frac{\varepsilon}{2})\geq \frac{\delta}{2}}\}|$
$$+ |{\{k\in [\alpha_{n},\beta_{n}]: P(|X_{k}-Y|\geq \frac{\varepsilon}{2})\geq \frac{\delta}{2}}\}|$$\\

$\Rightarrow (\beta_{n}-\alpha_{n})\leq |{\{k\leq [\beta_{n}]: P(|X_{k}-X|\geq \frac{\varepsilon}{2})\geq \frac{\delta}{2}}\}|+ |{\{k\in [\alpha_{n},\beta_{n}]: P(|X_{k}-Y|\geq \frac{\varepsilon}{2})\geq \frac{\delta}{2}}\}|$\\

$\Rightarrow (\frac{\beta_{n}-\alpha_{n}}{\beta_{n}})^{\gamma_1}\leq \frac{1}{[\beta_{n}]^{\gamma_1}}|{\{k\leq [\beta_n]: P(|X_{k}-X|\geq \frac{\varepsilon}{2})\geq \frac{\delta}{2}}\}|+ \frac{1}{\beta_{n}^{\gamma_1}}|{\{k\in [\alpha_{n},\beta_{n}]: P(|X_{k}-Y|\geq \frac{\varepsilon}{2})\geq \frac{\delta}{2}}\}|$\\

$\Rightarrow (\frac{\beta_{n}-\alpha_{n}}{\beta_{n}})^{\gamma_1}\leq \frac{1}{[\beta_{n}]^{\gamma_1}}|{\{k\leq [\beta_n]: P(|X_{k}-X|\geq \frac{\varepsilon}{2})\geq \frac{\delta}{2}}\}|$
$$+ (\frac{\beta_n-\alpha_n+1}{\beta_n})^{\gamma_2}\frac{1}{(\beta_n-\alpha_n+1)^{\gamma_2}}|{\{k\in [\alpha_n,\beta_n]: P(|X_{k}-Y|\geq \frac{\varepsilon}{2})\geq \frac{\delta}{2}}\}|$$\\

$\Rightarrow (1-\frac{\alpha_n}{\beta_n})^{\gamma_1}\leq \frac{1}{[\beta_{n}]^{\gamma_1}}|{\{k\leq [\beta_n]: P(|X_{k}-X|\geq \frac{\varepsilon}{2})\geq \frac{\delta}{2}}\}|$
$$+ (1-\frac{\alpha_n}{\beta_n}+\frac{1}{\beta_n})^{\gamma_2}\frac{1}{(\beta_n-\alpha_n+1)^{\gamma_2}}|{\{k\in [\alpha_n,\beta_n]: P(|X_{k}-Y|\geq \frac{\varepsilon}{2})\geq \frac{\delta}{2}}\}|.$$\\

Taking $n\rightarrow \infty$ on both sides we see that the left hand side does not tend to zero since $\liminf (\frac{\beta_n}{\alpha_n})>1$ but the right hand side tends to zero. This is a contradiction. So we must have $P\{X=Y\}=1.$\\

\section{\textbf{$\alpha\beta$-statistical convergence of order $\gamma$ in $r^{\mbox{th}}$ expectation}}

\noindent{\textbf{Definition 3.1.}} Let $(S,\triangle,P)$ be a probability space and $\{X_{n}\}_{n\in \mathbb{N}}$ be a sequence of random variables where each $X_{n}$ is defined on the same sample space $S$ (for each $n$) with respect to a given class of events $\triangle$ and a given probability function  $P:\triangle \rightarrow \mathbb{R}.$ Then the sequence $\{X_{n}\}_{n\in \mathbb{N}}$ is said to be $\alpha\beta$-statistically convergent of order $\gamma$ (where $0<\gamma\leq 1$) in $r^{\mbox{th}}$ expectation to a random variable $X$ (where $X:S\rightarrow \mathbb{R}$) if for any $\varepsilon>0$ $$ {\lim_{n\rightarrow \infty}}\frac{1}{(\beta_n-\alpha_n+1)^{\gamma}}|\{k\in [\alpha_n,\beta_n]: E(|X_{k}-X|^{r})\geq \varepsilon\}|=0,$$
  provided $E(|X_{n}|^{r})$ and $E(|X|^{r})$ exists $\forall$ $n\in \mathbb{N}.$ In this case we write $S_{\alpha\beta}^{\gamma}-\lim  E(|X_{n}-X|^{r})=0$ or by $X_{n}\xrightarrow{ES_{\alpha\beta}^{\gamma,r}} X.$  The class of all  sequences of random variables which are $\alpha\beta$-statistically convergent of order $\gamma$ in $r^{\mbox{th}}$ expectation is denoted simply by $ES_{\alpha\beta}^{\gamma,r}.$\\

\noindent{\textbf{Theorem 3.1.}} Let $X_{n}\xrightarrow{ES_{\alpha\beta}^{\gamma,r}} X $ (for any $r>0$ and $0<\gamma\leq 1$). Then $X_{n}\xrightarrow{PS_{\alpha\beta}^{\gamma}} X,$ i.e., $\alpha\beta$-statistical convergence of order $\gamma$ in $r^{\mbox{th}}$ expectation implies  $\alpha\beta$-statistical convergence of order $\gamma$ in probability.\\

\noindent{\textbf{Proof :}} The proof can be easily obtained by using Bienayme-Tchebycheff's inequality.\\

The following example shows that in general the converse of Theorem $3.1.$ is not true.\\

\noindent{\textbf{Example 3.1.}} We consider the sequence of random variables $\{X_{n}\}_{n\in \mathbb{N}}$ defined by,\\
\begin{equation*}
X_{n}\in
\begin{cases}
\{0,1\} ~~ \mbox{with probability} ~~ P(X_{n}=0)=P(X_{n}=1)~~ \mbox{if} ~~ n=m^2~~\mbox{for some}~~m\in \mathbb{N}\\
\{0,n\} ~~ \mbox{with probability} ~ P(X_{n}=0)=(1-\frac{1}{n^{r}})~~ \mbox{and} ~~ P(X_{n}=1)=\frac{1}{n^{r}},~~\mbox{if}~~ n\neq m^2~~\mbox{for any}~~m \in \mathbb{N}
\end{cases}
\end{equation*}
Where $r>0$. Now let $0<\varepsilon,\delta<1.~~ \mbox{Then we have}$
$$P(|X_{n}-0|\geq \varepsilon)=\frac{1}{2} ~~ \mbox{if} ~~ n=m^2~~\mbox{for some}~~m\in \mathbb{N}$$ and
$$P(|X_{n}-0|\geq \varepsilon)=\frac{1}{n^{r}},~~\mbox{if}~~n\neq m^2~~\mbox{for any}~~m \in \mathbb{N}.$$
Let $\gamma=\frac{1}{2}$, $\alpha_{n}=((n-1)^2+1)$, $\beta_{n}=n^2$. Then we have the inequality, $$\frac{1}{\sqrt{2n-1}}|\{k\in [(n-1)^2+1,n^2]: P(|X_{n}-0|\geq \varepsilon)\geq \delta\}|=(\frac{1}{\sqrt{2n-1}}+\frac{d}{\sqrt{2n-1}})\rightarrow 0 ~~\mbox{as}~~n\rightarrow \infty$$
where $d$ is a finite positive integer. So $X_{n}\xrightarrow {PS_{\alpha\beta}^{\frac{1}{2}}}0$.\\
%

But
\begin{eqnarray*}
E(|X_{n}-0|^{r})= & \frac{1}{2} & ~~ \mbox{if} ~~ n=m^2~~\mbox{for some}~~m\in \mathbb{N}\\
=& 1 &~~\mbox{if}~~ n\neq m^2~~\mbox{for any}~~m \in \mathbb{N}.
\end{eqnarray*}
This shows that $S_{\alpha\beta}^{\frac{1}{2}}$-$\lim E(|X_{n}-0|^{r})\neq 0$. i.e $X_{n}$ is not $\alpha\beta$-statistically convergent of order $\frac{1}{2}$ in $r^{\mbox{th}}$ expectation to $0$.\\

\noindent{\textbf{Theorem 3.2.}} Let $\{X_{n}\}_{n\in \mathbb{N}}$ be a sequence of random variables such that $P(|X_{n}|\leq M)=1$ for all $n$ and some constant $M>0.$ Suppose that $X_{n}\xrightarrow{PS_{\alpha\beta}^{\gamma}} X.$ Then $X_{n}\xrightarrow{ES_{\alpha\beta}^{\gamma,r}} X$ for any $r>0.$\\


 \noindent{\textbf{Theorem 3.3.}} (i) Let $X_{n}\xrightarrow{ES_{\alpha\beta}^{\gamma,r}} X$ and $X_{n}\xrightarrow{ES_{\alpha\beta}^{\gamma,r}} Y$ (for all $r>0$ and $0<\gamma\leq 1$). Then $P(X=Y)=1$ provided $(X-X_{n})\geq 0$ and $(Y_{n}-Y)\geq 0$.\\
(ii) Let $X_{n}\xrightarrow{ES_{\alpha\beta}^{\gamma,r}} X$ and $Y_{n}\xrightarrow{ES_{\alpha\beta}^{\gamma,r}} Y$ (for all $r>0$ and $0<\gamma\leq 1$). Then $(X_{n}+Y_{n})\xrightarrow{ES_{\alpha\beta}^{\gamma}} (X+Y)$ provided $(X-X_{n})\geq 0$ and $(Y_{n}-Y)\geq 0$.\\



\section{\textbf{$\alpha\beta$-statistical convergence of order $\gamma$ in distribution}}

\noindent{\textbf{Definition 4.1.}} Let $(S,\triangle,P)$ be a probability space and $\{X_{n}\}_{n\in \mathbb{N}}$ be a sequence of random variables where each $X_{n}$ is defined on the same sample space $S$ (for each $n$) with respect to a given class of events $\triangle$ and a given probability function  $P:\triangle \rightarrow \mathbb{R}.$ Let $F_{n}(x)$ is the distribution function of $X_{n}$ $\forall~ n\in \mathbb{N}$. If there exist a random variable $X$ whose distribution function is $F(x)$ such that the sequence $\{F_{n}(x)\}_{n\in \mathbb{N}}$ is $\alpha\beta$-statistically convergent of order $\gamma$ to $F(x)$ at every point of continuity $x$ of $F(x)$, then $\{X_{n}\}_{n\in \mathbb{N}}$ is said to be $\alpha\beta$-statistically convergent of order $\gamma$ in distribution to $X$ and we write $X_{n}\xrightarrow {DS_{\alpha\beta}^{\gamma}}X$.\\

\noindent{\textbf{Theorem 4.1.}} Let $\{X_{n}\}_{n\in \mathbb{N}}$ be a sequence of random variables. Also let $f_{n}(x)=P(X_{n}=x)$ be the probability mass function of $X_{n}~\forall~n \in \mathbb{N}$ and $f(x)=P(X=x)$ be the probability mass function of $X$. If $f_{n}(x)\xrightarrow {S_{\alpha\beta}^{\gamma}}f(x)~\forall~x$ then $X_{n}\xrightarrow {DS_{\alpha\beta}^{\gamma}}X$.\\

\noindent{\textbf{Proof :}} Proof is straightforward, so omitted.\\

\noindent{\textbf{Proposition 4.1.}} Let $\{a_{n}\}_{n\in \mathbb{N}}$, $\{b_{n}\}_{n\in \mathbb{N}}$ be two sequences of real numbers such that $a_{n}\leq b_{n}~ \forall~n \in \mathbb{N}$. Then $$S_{\alpha\beta}^{\gamma}-\underline{\lim} a_{n}\leq S_{\alpha\beta}^{\gamma}-\underline{\lim} b_{n}~\mbox{and}~S_{\alpha\beta}^{\gamma}-\overline{\lim} a_{n}\leq S_{\alpha\beta}^{\gamma}-\overline{\lim} b_{n}.$$
Here $S_{\alpha\beta}^{\gamma}-\underline{\lim}$ and $S_{\alpha\beta}^{\gamma}-\overline{\lim}~$ denotes $~\alpha\beta$-statistical limit inferior of order $\gamma$ and $\alpha\beta$-statistical limit superior of order $\gamma$ of the respective real sequences and here we use the same definition as in \cite{fo1} but here natural density is replaced by the $\alpha\beta$ density of order $\gamma.$\\

\noindent{\textbf{Proof :}} Proof is straightforward, so omitted.\\

\noindent{\textbf{Theorem 4.2.}} Let $\{X_{n}\}_{n\in \mathbb{N}}$ be a sequence of random variables. If $X_{n}\xrightarrow {PS_{\alpha\beta}^{\gamma}}X$ then $X_{n}\xrightarrow {DS_{\alpha\beta}^{\gamma}}X$. That is, $\alpha\beta$-statistical convergence of order $\gamma$ in probability implies $\alpha\beta$-statistical convergence of order $\gamma$ in distribution.\\

\noindent{\textbf{Proof.}} Let $F_{n}(x)$ and $F(x)$ be the probability distribution functions of $X_{n}$ and $X$ respectively. Let $x<y$. Now $$(X\leq x)\subseteq (X_{n}\leq y,X\leq x)+(X_{n}>y,X\leq x)$$
$$\Rightarrow (X\leq x)\subseteq (X_{n}\leq y)+(X_{n}>y,X\leq x)$$
$$\Rightarrow P(X\leq x)\leq P(X_{n}\leq y)+P(X_{n}>y,X\leq x)$$
$$\Rightarrow P(X\leq x)\leq P(X_{n}\leq y)+P(|X_{n}-X|>y-x)~(\mbox{as}(X_{n}>y,X\leq x)\subseteq (|X_{n}-X|>y-x))$$
$$\Rightarrow F(x)\leq F_{n}(y)+ P(|X_{n}-X|>y-x)$$
$$\Rightarrow S_{\alpha\beta}^{\gamma}-\underline{\lim}F(x)\leq S_{\alpha\beta}^{\gamma}-\underline{\lim}F_{n}(y)~(\mbox{as}~ S_{\alpha\beta}^{\gamma}-\lim P(|X_{n}-X|>y-x)=0, ~\mbox{since}~ X_{n}\xrightarrow {PS_{\alpha\beta}^{\gamma}}X)$$
$$\Rightarrow F(x)\leq S_{\alpha\beta}^{\gamma}-\underline{\lim}F_{n}(y)$$
similarly following same kind of steps and by taking $y<z$ we get, $$S_{\alpha\beta}^{\gamma}-\overline{\lim}F_{n}(y)\leq F(z)$$
Let y be a point of continuity of the function $F(x)$. Then $$\displaystyle{\lim_{x\rightarrow y-}}F(x)=\displaystyle{\lim_{z\rightarrow y+}}F(z)=F(y)$$ Hence we get $S_{\alpha\beta}^{\gamma}-\lim F_{n}(y)=F(y)$. Hence the result follows.\\

 Now we will show that the converse of the Theorem 4.2 is not necessarily true i.e $\alpha\beta$-statistical convergence of order $\gamma$ in distribution may not implies $\alpha\beta$-statistical convergence of order $\gamma$ in probability. For this we will construct an example as follows :\\

\noindent{\textbf{Example 4.1.}} Let $\alpha_{n}=((n-1)^2 +1)$ and $\beta_{n}=n^2$ and consider random variables $X,X_n$ (where $n$ is the first $[{h_{r}}^c]$ integers in the interval $I_{r}=[\alpha_{r},\beta_{r}]$ where $h_r=(\beta_r-\alpha_r+1)$ and $0<c<1$) having identical distribution. Let $\gamma$ be a real number in $(0,1]$ such that $\gamma>c.$ Let the spectrum of the two dimensional random variables $(X_{n},X)$ be $(-1,0),(-1,1),(1,0),(1,1)$ with probability $$P(X_{n}=1,X=1)=0=P(X_{n}=-1,X=0)$$  $$P(X_{n}=-1,X=1)=\frac{1}{2}=P(X_{n}=1,X=0).$$ Hence, the marginal distribution of $X_{n}$ is given by $X_{n}=i~(i=-1,1)$, with p.m.f, $f_{X_n}(-1)=f_{X_n}(1)=\frac{1}{2}$ and the marginal distribution of $X$ is given by $X=i~(i=0,1)$, with p.m.f, $f_{X}(0)=f_{X}(1)=\frac{1}{2}.$\\

Next, we consider random variables $X,X_n$ (where $n$ is other than the first $[{h_{r}}^c]$ integers in the interval $I_{r}=[\alpha_{r},\beta_{r}]$) having identical distribution and the spectrum of the two dimensional random variables $(X_{n},X)$ be $(0,0),(0,1),(1,0),(1,1)$ with probability $$P(X_{n}=0,X=0)=0=P(X_{n}=1,X=1)$$  $$P(X_{n}=1,X=0)=\frac{1}{2}=P(X_{n}=0,X=1)$$ Hence, the marginal distribution of $X_{n}$ is given by $X_{n}=i~(i=0,1)$, with p.m.f, $f_{X_n}(0)=f_{X_n}(1)=\frac{1}{2}$ and the marginal distribution of $X$ is given by $X=i~(i=0,1)$, with p.m.f, $f_{X}(0)=f_{X}(1)=\frac{1}{2}.$ \\

Let $n$ is the first $[{h_{r}}^c]$ integers in the interval $I_{r}$ and $F_{n}(x)$ be the probability distribution function of $X_{n}$ then,
\begin{equation*}
F_{n}(x)=
\begin{cases}
0 ~~\mbox{if},~ x<-1,\\
\frac{1}{2}~~  \mbox{if},~ -1\leq x<1,~&\\
1~~  \mbox{if},~x\geq 1
\end{cases}
\end{equation*}

Next let, $n$ is other than the first $[{h_{r}}^c]$ integers in the interval $I_{r}$ and $F_{n}(x)$ and $F(x)$ be the probability distribution function of $X_{n}$ and $X$ respectively, then,
\begin{equation*}
F_{n}(x)=F(x)=
\begin{cases}
0 ~~\mbox{if},~ x<0,\\
\frac{1}{2}~~  \mbox{if},~ 0\leq x<1,~&\\
1~~  \mbox{if},~x\geq 1
\end{cases}
\end{equation*}

We consider the interval $[-1,0).$ It is sufficient to prove that the sequence $\{y_n\}_{n\in \mathbb{N}}$ define below is $\alpha\beta$-statistically convergent of order $\gamma$ to $0.$  Now we define a sequence $\{y_n\}_{n\in \mathbb{N}}$ by
\begin{equation*}
y_n=
\begin{cases}
\frac{1}{2} ~~ \mbox{if} ~~n ~\mbox{is the first}~[{h_{r}}^c]~\mbox{integers in the interval}~I_{r},&\\
0,~~\mbox{if}~~ n ~\mbox{is other than the first}~[{h_{r}}^c]~\mbox{integers in the interval}~I_{r}.
\end{cases}
\end{equation*}
It is quite clearly that $\{y_n\}_{n\in \mathbb{N}}$ is $\alpha\beta$-statistically convergent of order $\gamma$ to $0,$ this implies $F_{n}(x)\xrightarrow {S_{\alpha\beta}^\gamma}F(x)~ \forall~ x\in \mathbb{R}$ (but $\displaystyle{\lim_{n\rightarrow \infty}}F_n(x)\neq F(x)~\forall~x\in [-1,0)$ in ordinary sence). This shows that $X_{n}\xrightarrow {DS_{\alpha\beta}^{\gamma}}X$.\\

For any $0<\epsilon<1,$
\begin{equation*}
P(|X_{n}-X|\geq \epsilon)=
\begin{cases}
1 ~~ \mbox{if} ~~n ~\mbox{is the first}~~[{h_{r}}^c]~\mbox{integers in the interval}~I_{r},&\\
1,~~\mbox{if}~~ n ~\mbox{is other than the first}~~[{h_{r}}^c]~\mbox{integers in the interval}~I_{r}.
\end{cases}
\end{equation*}

 But this shows that ${S_{\alpha\beta}^\gamma}-\lim P(|X_{n}-X|\geq \epsilon)\neq 0$. This shows that the sequence $\{X_{n}\}_{n\in \mathbb{N}}$ is not $\alpha\beta$-statistically convergent of order $\gamma$ in probability to $X$.\\

\noindent{\textbf{Note 4.1.}} The statement of Theorem 4.2 does not depend on the limit infimum of the sequence  $q_{r}=\frac{\beta_{r}}{\alpha_{r}}$ of the  sequence $\alpha_{n}$ and $\beta_{n}.$ Even if $\liminf q_{r}=1$, Theorem 4.2 will hold good.\\

\noindent{\textbf{Theorem 4.3.}} Let $\{\alpha_n\}_{n\in \mathbb{N}}$ and $\{\beta_n\}_{n\in \mathbb{N}}$ be two increasing sequences of positive real numbers such that $\alpha_n\leq \beta_n\leq \alpha_{n+1}\leq \beta_{n+1}$~, $(\beta_n-\alpha_n)\rightarrow \infty$ as $n\rightarrow \infty$ and $0<\gamma_1\leq \gamma_2\leq 1$. If $\liminf (\frac{\beta_n}{\alpha_n})>1$ then $X_{n}\xrightarrow {PS^{\gamma_1}}X$ implies $X_{n}\xrightarrow {DS_{\alpha\beta}^{\gamma_2}}X. $\\

\noindent{\textbf{Proof :}} Proof is straightforward, so omitted.\\







\end{document}